\documentclass[a4paper]{article}
\usepackage[margin=2cm]{geometry}

\usepackage{amsmath}
\usepackage{amssymb}
\usepackage{amsfonts}
\usepackage{enumitem}
\usepackage{nicefrac}
\usepackage{cite}
\usepackage{amsthm}
\usepackage{hyperref,url}

\usepackage[ruled,vlined,linesnumbered,algosection,algo2e]{algorithm2e} 
\usepackage[capitalize,noabbrev,nameinlink]{cleveref} 

\usepackage{xcolor}
\definecolor{cb2blue}{RGB}{55,126,184}
\definecolor{cb2green}{RGB}{77,175,74}
\definecolor{cb2red}{RGB}{228,26,28}

\theoremstyle{plain}
\newtheorem{theorem}{Theorem}[section]

\hypersetup{colorlinks,
	citecolor=cb2green,
	linkcolor=cb2red}
	
\newcommand{\emailLink}[1]{\textsc{email} \href{mailto:#1}{#1}}
\newcommand{\orcidLink}[1]{\textsc{orcid} \href{https://orcid.org/#1}{#1}}
\newcommand{\keywords}[1]{\par\noindent{\def\and{\unskip,\ }{\bf Keywords: }#1}\par}

\newcommand{\stepref}[1]{\hyperref[#1]{Step~\ref{#1}}}
\SetKw{KwOr}{or}
\SetKw{KwAnd}{and}

\SetCommentSty{mycommfont}
\makeatletter
	\newcommand{\algcolor}[2]{\colorbox{#1}{#2}}
	\newcommand{\algcoloremph}[1]{\algcolor{cb2blue!10!white}{#1}}
\makeatother

\usepackage{pgfplots}
\pgfplotsset{compat=newest}
\usepackage{tikz}
\usetikzlibrary{calc, intersections}
\usetikzlibrary{shapes.geometric}

\newcommand{\coloneqq}{:=}
\DeclareMathOperator*{\minimize}{minimize}

\DeclareMathOperator{\stt}{subject~to}
\newcommand{\func}[3]{#1 \colon #2 \to #3}
\newcommand{\ffunc}[3]{#1 \colon #2 \rightrightarrows #3}
\newcommand{\N}{\mathbb{N}}
\newcommand{\R}{\mathbb{R}}

\newcommand{\innprod}[2]{\langle #1, #2 \rangle}
\newcommand{\ang}[1]{\langle #1 \rangle}
\DeclareMathOperator{\projection}{\Pi}

\DeclareMathOperator{\indicator}{\iota}

\newcommand{\XX}{\mathcal{X}}
\newcommand{\limnormalcone}{\mathcal{N}^{\text{lim}}}
\DeclareMathOperator{\UB}{UB}
\DeclareMathOperator{\LB}{LB}

\usepackage{mathrsfs}
\newcommand{\resid}{\mathscr{R}}

\newcommand{\TheTitle}{Cutting plane methods with gradient-based heuristics}

\newcommand{\TheAuthorHB}{H\`oa T. B\`ui}
\newcommand{\TheAuthorADM}{Alberto De~Marchi}
\newcommand{\TheAffiliationADM}{%
	Institute of Applied Mathematics and Scientific Computing,
	Department of Aerospace Engineering,
	University of the Bundeswehr Munich, Germany%
}

\newcommand{\TheAffiliationCCODS}{%
	Curtin Centre for Optimisation and Decision Science,
	Curtin University,
	Australia%
}
\newcommand{\TheKeywords}{%
    Integer programming,
    outer approximation,
    cutting planes,
    continuous optimization,
    projected gradient algorithms%
}
\newcommand{\TheAcknowledgements}{This research was partially conducted during ADM’s visit to Curtin University, Perth, WA. The author gratefully acknowledges the support and hospitality of Ryan Loxton and the Centre for Optimisation and Decision Science.}

\begin{document}

\title{\bfseries \TheTitle}
\title{\bfseries \TheTitle}
\author{\TheAuthorHB\thanks{\TheAffiliationCCODS. \emailLink{hoa.bui@curtin.edu.au}, \orcidLink{0000-0002-1698-6383}.}%
	\and%
	\TheAuthorADM\thanks{\TheAffiliationADM. \emailLink{alberto.demarchi@unibw.de}, \orcidLink{0000-0002-3545-6898}.}%
}
\date{}

\maketitle

\begin{abstract}
    Cutting plane methods, particularly outer approximation, are a well-established approach for solving nonlinear discrete optimization problems without relaxing the integrality of decision variables. While powerful in theory, their computational performance can be highly variable. Recent research has shown that constructing cutting planes at the projection of infeasible points onto the feasible set can significantly improve the performance of cutting plane approaches.
    Motivated by this, we examine whether constructing cuts at feasible points closer to the optimal solution set could further enhance the effectiveness of cutting plane methods.
    We propose a hybrid method that combines the global convergence guarantees of cutting plane methods with the local exploration capabilities of first-order optimization techniques. Specifically, we use projected gradient methods as a heuristic to identify promising regions of the solution space and generate tighter, more informative cuts. We focus on binary optimization problems with convex differentiable objective functions, where projection operations can be efficiently computed via mixed-integer linear programming. By constructing cuts at points closer to the optimal solution set and eliminating non-optimal regions, the algorithm achieves better approximation of the feasible region and faster convergence. Numerical experiments confirm that our approach improves both the quality of the solution and computational efficiency across different solver configurations. This framework provides a flexible foundation for further extensions to more general discrete domains and offers a promising heuristic to the toolkit for nonlinear discrete optimization.
    
    \medskip
    
\end{abstract}

\section{Introduction}
The cutting plane method (a.k.a. outer approximation) \cite{kelley1960cutting,grossmann1986outer,westerlund1995extended} is one of the primary approaches for solving nonlinear discrete optimization problems without requiring continuous relaxation of the discrete variable set. 
The methodology has been extensively investigated also for mixed-integer problems \cite{fletcher1994solving,bonami2008algorithmic,belotti2013mixed}.
However, the computational performance of cutting plane methods exhibits significant variability and can deteriorate to brute-force efficiency levels in worst-case scenarios (see \cite{hijazi2014outer} for an example on the worst case scenarios of outer
approximation). This is because the convergence of cutting plane algorithms relies on the added cutting planes being sufficient to approximate the nonlinear components.
Since the seminal work of Grossmann and Kocis (1986) \cite{grossmann1986outer}, substantial efforts have concentrated on enhancing the strength and effectiveness of cutting planes.
A fundamental question in advancing cutting plane methodologies concerns where to best construct cutting planes. 
For instance, the extended cutting plane method \cite{westerlund1995extended} solves iteratively a MILP cutting plane model to construct one cutting plane at the optimal solutions of the MILP after each iteration. While each iteration requires solving an MILP to optimality, the resulting cuts are empirically tighter than those produced by an outer-approximation branch-and-cut framework \cite{bui2024solving}.
Moreover, the method typically adds fewer cuts, and the MILP-derived cuts reveal directions along which the current approximation is the most insufficient.
The extended supporting hyperplane \cite{kronqvist2016extended} follows the same fashion, but builds the initial cutting plane oracle by solving the continuous relaxation. Recent research \cite{westerlund2022projected} shows that constructing cuts as projected supporting hyperplane, projecting infeasible solutions from continuous relaxation to the feasible set, yields much tighter relaxation.
This technique involves generating multiple cuts for each solution of the mixed-integer nonlinear programming (MINLP) subproblems and projecting infeasible points toward the boundary of the feasible region.
Through this projection step, the feasible region is more accurately approximated, leading to improved algorithmic performance.
Following the same idea, we explore two potential avenues to improve the tightness of cutting planes: (1) cuts constructed at points that are closer to the optimal solutions set yield tighter approximations, and (2) cuts that are designed to eliminate non-optimal solutions are more effective in reducing the feasible region.
In this paper, we use gradient-based heuristics to find tighter cutting planes that can improve the characterization of the optimal solution set.

We consider the following nonlinear discrete minimization problem
\begin{equation}\label{eq:P}\tag{NDP}
	\minimize_{x \in \XX}
	f(x),
\end{equation}
where $\XX \coloneqq \mathbb{Z}^n \cap P$ is the feasible set,
$f$ is a convex differentiable function,
and $P\subset\R^n$ is a nonempty bounded polyhedral set.

Gradient-based methods are recently getting attentions as possible solvers for \eqref{eq:P}, dealing directly with the discrete set $\XX$ as a nonconvex set (e.g., \cite{sotirov2020solving,ghaffari2024convergence,bragin2018scalable}).
They are attractive from a scalability standpoint: iterations involve only gradient and projection operations, and the number of iterations is generally independent of the problem size \cite{cartis2011evaluation}. 
However, first-order methods are local in their nature:
they can quickly identify local solutions but often hinders convergence to global optimal solutions.
Conversely, branch and bound and outer approximation techniques guarantee global optimality but they may need many iterations to effectively reduce the optimality gap.
In this paper, we combine the two approaches to alleviate these issues.
For the sake of clarity, we will illustrate our strategy by focusing on the binary case and by consolidating projected gradient (first-order) and cutting planes (global) methods, but other settings and combinations are possible.

\subsection{Cutting plane methods and motivations}

Given a set of point $\mathcal{A}\subset \R^n$, the cutting plane model \ref{eq:thetaProblem} is defined as:
\begin{align}
	\minimize_{x\in\XX ,\, \theta\in\R}\quad{}& \theta \label{eq:thetaProblem}\tag{CP($\mathcal{A}$)}\\
	\stt\quad{}& \theta \geq f(x_a) + \innprod{\nabla f(x_a)}{x-x_a} \qquad\forall x_a\in\mathcal{A}. \nonumber
\end{align}
The set $\mathcal{A}$ is called cutting plane oracle.
For each $x_a\in\mathcal{A}$, the constraint $\theta \geq f(x_a) + \innprod{\nabla f(x_a)}{x-x_a}$ (called cutting plane) is the tangent plane of $f$ constructed at $x_a \in \mathcal{A}$.
When the function $f$ is convex, the cutting plane model \ref{eq:thetaProblem} provides an upper bound for problem \eqref{eq:P}.
Additionally, when $\mathcal{A} = \XX$, the two problems \eqref{eq:P} and \ref{eq:thetaProblem} are equivalent in the sense that $x^*$ is an optimal solution of the nonlinear problem \eqref{eq:P} if and only if $(x^*, f(x^*))$ is an optimal solution of the cutting plane model \ref{eq:thetaProblem}.
The cutting plane method (\cref{alg:cp}) often starts with an initial set $\mathcal{A}$, the cutting plane model \ref{eq:thetaProblem} is solved to optimality to obtain a solution $(x_a, \theta_a)$, then a new cutting plane is constructed at $x_a$ with $\theta_a$ providing a valid upper bound for the nonlinear model \eqref{eq:P}.
The algorithm terminates when there is no new cut generated.
The convergence of cutting plane method to optimal solutions relies on two criteria:
(1) every added cut is valid (i.e., it does not exclude the optimal solution $(x^*, f(x^*))$, and
(2) the added cut must remove at least one non-optimal point \cite{bui2025cutting}. 

\begin{algorithm2e}
	\DontPrintSemicolon%
	\KwData{$x^0 \in \XX, \varepsilon > 0$}
	\KwResult{optimal solution $\bar{x} \in \XX$}
    
    $\mathcal{A}_0 \gets \{x^0\}$, $\UB_k \gets f(x^0)$ and $\LB_k \gets - \infty$\;

    \For{$k \gets 0$ \KwTo $\infty$}{%
        \lIf{$\UB_k - \LB_k \leq \varepsilon$}{\Return $x^k$\tcp*[f]{termination}}
        
        $x^{k+1}, \theta_{k+1} \gets$ \texttt{MILP\_solver}(\hyperref[eq:thetaProblem]{CP($\mathcal{A}_k$)})
		\label{step:cp:lowerBoundProblem}\tcp*{solve cutting plane model}
        
        $\mathcal{A}_{k+1} \gets \mathcal{A}_k\cup \{x^k\}$
        \label{step:cp:update_cuts}\tcp*{add cutting plane}
            
        $\UB_{k+1} \gets \min\{\UB_k, f(x^{k+1})\}$ and $\LB_{k+1} \gets \theta_{k+1}$\label{step:cp:update_bounds}\tcp*{update bounds}
    }	
	\caption{Standard cutting-plane method (CPM) for \eqref{eq:P}.}%
	\label{alg:cp}%
\end{algorithm2e}

The key step in \cref{alg:cp} is updating the next iterate and its associated cutting plane. In the standard cutting-plane method, the cut is generated at the optimal solution $x^{k+1}$ of \hyperref[eq:thetaProblem]{CP($\mathcal{A}_k$)}.
Although $x^{k+1}$ solves the cutting-plane model, its objective value $f(x^{k+1})$ is not necessarily smaller than that of other nearby solutions.
To tighten the upper bound and construct a better cut, we invoke a local solver to obtain an improved feasible solution while preserving the search direction obtained from the cutting plane model.

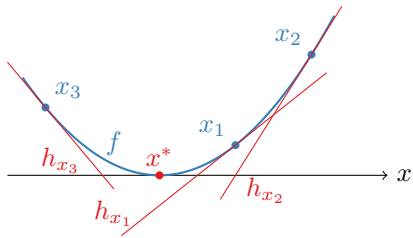
\begin{figure}
    \centering%
    \begin{tikzpicture}
	    \draw[->] (-2,0) -- (3,0) node[right] {$x$};
	    
	    \draw[cb2blue, thick, smooth] plot[domain=-1.8:2.3, samples=100] (\x, {0.4*\x*\x});
	    
	    \coordinate (P) at (1, 0.4);
	    \fill[cb2blue] (P) circle (1.5pt);
	    \draw[cb2red] plot[domain=-0.5:2.2] (\x, {0.8*(\x-1) + 0.4});
	    \node[cb2blue, above left] at (P) {$x_1$};
	    
	    \coordinate (L) at (2, 1.6);
	    \fill[cb2blue] (L) circle (1.5pt);
	    \draw[cb2red] plot[domain=0.8:2.4] (\x, {1.6*(\x-2) + 1.6});
	    \node[cb2blue, above left] at (L) {$x_2$};
	    
	    \coordinate (L) at (-1.5, 0.9);
	    \fill[cb2blue] (L) circle (1.5pt);
	    \draw[cb2red] plot[domain=-2:-0.6] (\x, {-1.2*(\x+1.5) + 0.9});
	    \node[cb2blue, above right] at (L) {$x_3$};
	    
	    \node[align=center, cb2red] at (-0.6, -0.5) {$h_{x_1}$};
	    \node[align=center, cb2red] at (1.4, -0.2) {$h_{x_2}$};
	    \node[align=center, cb2red] at (-1.3, 0.2) {$h_{x_3}$};
	     \coordinate (Q) at (0, 0);
	    \fill[cb2red] (Q) circle (1.5pt);
	    \node[cb2red, above] at (Q) {$x^*$};
	    \node[cb2blue, above] at (-0.6,0.144) {$f$};
	\end{tikzpicture}
    \caption{The tangent plane constructed at $x_1$ (points that are close to the optimal solution $x^*$) provides tighter bounds than the tangent plane constructed at $x_2$. Combining cuts at points $x_1$ and $x_3$, on opposite sides of $x^*$, also generates a tight relaxation.}%
    \label{fig:cvx_illustration}%
\end{figure}

Given an upper bound $\UB$ obtained in \stepref{step:cp:update_bounds} of \cref{alg:cp},
we invoke a local solver to address the problem
\begin{align}\label{eq:local}
    \minimize_{x\in \mathcal{C}_{\mathcal{A}, \text{UB}}}&\quad f(x),
\end{align}
where
\begin{subequations}
    \label{eq:cutsDefinition}
    \begin{align}
        \mathcal{C}_{\mathcal{A},\UB}
        \coloneqq{}&
        \bigcap_{x_a\in \mathcal{A}} F_{x_a, \UB}, \\
        F_{x_a, \UB}
        \coloneqq{}&
        \left\{
    	x\in\R^n
    	\,\middle\vert\,
    	f(x_a) + \innprod{\nabla f(x_a)}{x-x_a} \leq \UB - \tau
    	\right\},
    \end{align}
\end{subequations}
and $\tau \geq 0$ is an acceptable tolerance.
To simplify the notation, where $\mathcal{A}$ and $\UB$ are clearly stated, we can write $\mathcal{C}$ instead of $\mathcal{C}_{\mathcal{A},\UB}$.
If $f(x_a) > \UB$ then the solution $(x_a, f(x_a))$ is excluded by the cutting plane. 
Combining local searches with the cutting plane method \cref{alg:cp}, we propose \cref{alg:BUI}, a gradient-based heuristic scheme for \eqref{eq:P} that incorporates local steps to strengthen cut quality.

\begin{algorithm2e}
	\DontPrintSemicolon%
	\KwData{$x^0 \in \XX$, $\varepsilon>0$, \algcoloremph{\texttt{use\_local\_solver}, \texttt{use\_offset}, \texttt{use\_LB\_cuts}}\;}
	\KwResult{$\varepsilon$-optimal $\bar{x} \in A$ for \eqref{eq:P}\;}
    \algcoloremph{choose $\tau_0>0$, $\kappa_\tau\in(0,1)$ and $\kappa_g\in(0,1]$}\;
	$x^0_{\UB} \gets x^0$, $\UB_0 \gets f(x^0)$, $\mathcal{A}_0 \gets \{x^0\}$, \algcoloremph{$\texttt{increase\_offset}\gets\texttt{use\_offset}$}\;
	\For{$k \gets 0$ \KwTo $\infty$}{
		$x^k_{\LB}, \LB_k \gets$ \texttt{MILP\_solver}(\hyperref[eq:thetaProblem]{CP($\mathcal{A}_k$)})
		\label{step:BUI:lowerBoundProblem}\tcp*{lower bound problem}
		\lIf{$\UB_k - \LB_k \leq \varepsilon$\label{step:BUI:terminationCondition}}{\Return $x^k_{\UB}$\label{step:BUI:termination}\tcp*[f]{termination}}
		\eIf{\texttt{use\_local\_solver}}{%
            \eIf{\texttt{use\_offset}}{%
                \algcoloremph{$\tau_k \gets \min\{\tau_k, \kappa_g(\UB_k - \LB_k)\}$}\label{step:BUI:tauSaturation}\;
                \While{\texttt{true}}{%
                    \algcoloremph{$\mathcal{C}_k\gets$ cuts \eqref{eq:cutsDefinition} at $\mathcal{A}_k$ with bound $\UB_k$ and offset $\tau_k$}\; 
                    \algcoloremph{\lIf{$\XX\cap\mathcal{C}_k \neq \emptyset$}{\textbf{break}}}\tcp*{accept offset}
                    \algcoloremph{$\tau_k \gets \kappa_\tau \tau_k$ and $\texttt{increase\_offset}\gets\texttt{false}$}\tcp*{reduce offset}
                }
            }{%
                \algcoloremph{$\mathcal{C}_k\gets$ cuts \eqref{eq:cutsDefinition} at $\mathcal{A}_k$ with bound $\UB_k$}\;
            }%
            \algcoloremph{$x^{k+1} \gets \texttt{local\_solver}(f, \XX\cap\mathcal{C}_k)$ starting from $x^k_{\LB}$}\tcp*{local method}\label{step:BUI:localStep}
                \algcoloremph{\lIf{\texttt{increase\_offset}}{$\tau_{k+1}\gets \tau_k/\kappa_\tau$ \textbf{else} $\tau_{k+1}\gets\tau_k$}}
        }{%
            $x^{k+1} \gets x^k_{\LB}$
        }
		\If{$f(x^{k+1}) \leq \UB_k$}{%
			$x^{k+1}_{\UB} \gets x^{k+1}$ and $\UB_{k+1} \gets f(x^{k+1})$\tcp*{update best value}
		}
		\Else{%
			$x^{k+1}_{\UB} \gets x^k_{\UB}$ and $\UB_{k+1} \gets \UB_k$\;
		}
		$\mathcal{A}_{k+1}\gets \mathcal{A}_k\cup \{x^{k+1}\}$\tcp*{add cut}
		\If{\texttt{use\_LB\_cuts} $\wedge$ \algcoloremph{$\innprod{\nabla f(x^{k+1})}{x^k_{\LB}-x^{k+1}}\leq0$}}{
				\algcoloremph{$\mathcal{A}_{k+1}\gets \mathcal{A}_{k+1}\cup \{x^k_{\LB}\}$}\tcp*{add LB cut}
		}
	}
	\caption{Proposed cutting plane scheme with gradient-based heuristics for \eqref{eq:P}. Without heuristics (shaded steps) it recovers the classical cutting plane scheme in \cref{alg:cp}.}%
	\label{alg:BUI}%
\end{algorithm2e}

\paragraph{Outline}
The remainder of the paper is organized as follows.
\Cref{sec:heuristics} suggests first-order methods for generating tighter local cuts while preserving the global information provided by the MILP solve, thereby improving the efficiency of cutting-plane methods.
We place particular emphasis on the projected gradient method as a computationally inexpensive local solver for \eqref{eq:local} and analyze its convergence properties.
\Cref{sec:cut_tightenings} introduces additional cut-strengthening techniques; in particular, when the local solver’s solution deviates from the MILP solution, we add an additional cut obtained from the MILP solution.
\Cref{sec:numerical_experiments} presents numerical experiments that quantify the contribution of each local step to cut strengthening and overall performance.

\section{Heuristics: Local first-order methods on binary domain}\label{sec:heuristics}

Given a incumbent point $x^k$, it is reasonable to greedily seek a better one, for instance using ideas from continuous optimization to run an informed ``local'' search.
With this goal in mind, we consider auxiliary problems of the form
\begin{equation}\label{eq:localProblem}
	\minimize_{x\in A}\quad f(x)
\end{equation}
where function $\func{f}{A}{\R}$ is continuously differentiable and set $A\subset \R^n$ is possibly nonconvex but easy to project onto.
In the context of our CPM, the constraint set $A$ will be the intersection of $\XX$ with the current cuts $\mathcal{C}_k$ from \eqref{eq:cutsDefinition}; additional tightenings can be included too.
Several numerical methods from (nonsmooth and/or nonconvex) continuous optimization can provide convergence guarantees,
but typically only to ``local'' minimizers (or merely stationary points).
Well-known strategies are projected gradient descent, trust-region and Frank-Wolfe schemes, which rely on similar oracles.
We focus here on the projected gradient approach because of its simplicity (and practical relevance in the case of binary variables).
Frank-Wolfe schemes build upon a relatively practical linear minimization oracle but often some form of convexity in the problem plays a role in establishing convergence \cite{hendrych2025convex}.
In contrast, trust-region methods remain practical for a broader class of problems (mixed-integer variables, for instance) and offer greater versatility \cite{demarchi2025mixed}; these are subject of ongoing research, see \cref{sec:trust_region}.

\subsection{Projected gradient steps}\label{sec:projected_gradient}

The projected gradient method (PGM) applied to \eqref{eq:localProblem} starts from some $x^0\in A$ and employs an iterative procedure of the form
\begin{equation}\label{eq:PGiteration}
	x^{j+1} \in \projection_A\left(x^j - \gamma_j \nabla f(x^j) \right),
    \quad j=0,1,\ldots
\end{equation}
with some suitable stepsize $\gamma_j>0$.
The projection operator $\ffunc{\projection_A}{\R^n}{A}$ involves the squared Euclidean norm but it boils down to a MILP thanks to the binary structure of $A\subseteq \{0,1\}^n$:
\begin{align*}
	\projection_A(z)
	\coloneqq
	\arg\min_{x\in A} ~ \|x-z\|^2
	=
	\arg\min_{x\in A} ~ \innprod{1-2z}{x}
    .
\end{align*}
Therefore, we can assume its evaluation to be cheap enough (relative to that of $\nabla f$), although still NP-hard in general.
Note that the update rule \eqref{eq:PGiteration} corresponds to minimizing a quadratic model of $f$ around $x^j$ over the feasible set $A$:
\begin{equation*}
	x^{j+1} \in \arg\min_{x\in A}\left( f(x^j) + \innprod{\nabla f(x^j)}{x-x^j} + \frac{1}{2\gamma_j}\|x-x^j\|^2  \right) .
\end{equation*}
The stepsize $\gamma_j$ can be chosen by backtracking until either
\begin{equation}
	\label{eq:PGsuffdecrease}
	f(x^{j+1})
	\leq
	f(x^j) - \alpha \frac{\| x^{j+1}-x^j\|^2}{2\gamma_j}
\end{equation}
or
\begin{equation}
	\label{eq:PGqualityratio}
	f(x^{j+1})
	\leq
	f(x^j) - \alpha \innprod{\nabla f(x^j)}{x^j-x^{j+1}}
\end{equation}
holds, for some prescribed $\alpha\in(0,1)$, see \cref{alg:PGM}.
Both \eqref{eq:PGsuffdecrease} and \eqref{eq:PGqualityratio} encode the requirement that sufficient decrease in the cost $f$ is attained when updating from $x^j$ to $x^{j+1}$.
Common mechanisms to speed up the backtracking routine are the so-called \emph{spectral stepsize} and \emph{nonmonotone linesearch}.
The former provide a principled value for $\gamma_j$ based on local curvature estimates,
whereas the latter relaxes the requirement of monotone decrease for the cost $f$.
These tricks aim to reduce the number of times both \eqref{eq:PGsuffdecrease} and \eqref{eq:PGqualityratio} are violated;
see \cite[Alg. 3.1]{demarchi2023monotony} for more details.

\begin{algorithm2e}
	\DontPrintSemicolon%
	\KwData{$x^0 \in A$}
	\KwResult{critical $\bar{x} \in A$}
    Choose $\alpha, \beta \in (0,1)$ and $\gamma_0 \in (0,\infty)$\;
	\For{$j \leftarrow 0$ \KwTo $\infty$}{
        Set $\gamma_{j+1}\gets \gamma_j$\;
		\While{true}{
			Compute $x^{j+1} \in \projection_A(x^j - \gamma_{j+1} \nabla f(x^j))$\tcp*{subproblem}\label{step:PGM:projection}
            \lIf{$x^{j+1}=x^j$}{\Return $\bar{x} \gets x^j$\tcp*[f]{termination}}
			\lIf{\eqref{eq:PGsuffdecrease} or \eqref{eq:PGqualityratio}\label{step:PGM:linesearch:if}}{\textbf{break}\tcp*[f]{sufficient decrease}}
			Set $\gamma_{j+1} \gets \beta \gamma_{j+1}$\label{step:PGM:linesearch:backtrack}\tcp*{backtracking}
		}
	}
	\caption{Projected gradient method for \eqref{eq:localProblem}}%
	\label{alg:PGM}%
\end{algorithm2e}

\subsection{Convergence}

The convergence theory for \cref{alg:PGM} patterns that of proximal gradient methods, e.g. \cite{themelis2018forward,demarchi2022proximal,demarchi2023monotony},
but stronger properties can be established thanks to the additional structure of set $A$ in \eqref{eq:localProblem}.
Let us consider first an optimality concept for \eqref{eq:localProblem}:
a point $\bar{x}\in \R^n$ is called \emph{critical} for \eqref{eq:localProblem} if
\begin{equation}\label{eq:criticality}
	\exists \eta>0\colon~ \bar{x} \in \projection_A(\bar{x} - \eta \nabla f(\bar{x})) .
\end{equation}
Note that feasibility is necessary for criticality, namely $\bar{x}\in A$ must hold too.
See \cite[\S 3]{themelis2018forward} and \cref{sec:local_optimality} for more details.

Starting from a feasible point $x^0 \in A$,
\cref{alg:PGM} cannot produce an estimate $\bar{x}\in A$ worse than $x^0$ for \eqref{eq:localProblem},
since it generates feasible iterates while monotonically decreasing the objective value.
Convergence guarantees to stationary points follow, e.g., from \cite[\S 4.4]{demarchi2022proximal}, but \cref{alg:PGM} reveals peculiar properties when applied to \eqref{eq:localProblem}:
(i) accumulation points are not only stationary but also critical,
(ii) its termination is exact and finite, for arbitrary initial points.
Both features are due to discreteness of $A \subseteq \XX$ in our setting;
proof in \cref{sec:local_optimality}.
\begin{theorem}\label{prop:PGM:finiteExactTermination}
	Consider problem \eqref{eq:localProblem} with function $\func{f}{A}{\R}$ continuously differentiable and set $A \subset \N^n$ nonempty and compact.
	Then, starting from any $x^0 \in A$, \cref{alg:PGM} generates a finite sequence $\{x^j\}_{j=0,\ldots,J} \subseteq A$ and terminates with a critical point $\bar{x} \in A$ for \eqref{eq:localProblem} with $f(\bar{x}) \leq f(x^0)$.
\end{theorem}

\section{Additional cuts tightening}\label{sec:cut_tightenings}

The optional call to the local solver at \stepref{step:BUI:localStep} of \cref{alg:BUI} is intended to improve the objective value while maintaining feasibility.
Seeking to boost quality of such a local search, we tighten the feasible set $A$ with additional constraints, other than cuts $\mathcal{C}_k$.
As it is common with heuristics, there is a trade-off at play here:
the constraint set $A$ for the local search can be tightened to (hopefully) yield better values, but at the cost of a likely more computationally expensive projection operator.
Two techniques are presented in the following \cref{sec:cut_offset,sec:cut_lower_bound}, which aim to pass additional ``global'' information to help the solver eluding ``local'' solutions.
These optional and non-exclusive features will then be compared in the numerical section.

\subsection{Offset: trade-off and search}\label{sec:cut_offset}

One simple approach to overcome the incumbent solution $x_{\UB}^k$ with \eqref{eq:localProblem} is to impose, along with the constraint $x\in A$, that
$f(x) \leq \UB_k - \tau$
holds for some prescribed $\tau>0$.
However, this constraint is nonlinear, and expensive to handle.
Therefore, we approximate the function $f$ by a set of linear functions, exploiting convexity of $f$. 
In the spirit of CPM, we refine the approximation of $f$ using the set $\mathcal{A}_k$ and thus require the more practical constraints
\begin{align}
	\label{eq:offset_constraints}
    \forall x_a \in \mathcal{A}_k \colon\quad
    \ang{\nabla f(x_a), x}
	\leq
	\UB_k - \tau - f(x_a) + \ang{\nabla f(x_a), x_a}
\end{align}
along with $x\in \XX$.
Conditions \eqref{eq:offset_constraints}, already mentioned in \eqref{eq:cutsDefinition}, complement the original constraint $x\in \XX$ to eliminate all feasible solutions that do not improve the objective value.

For large values of $\tau$ there may be no $x\in\XX$ that satisfies \eqref{eq:offset_constraints},
while small values of $\tau$ may yield little improvement over $\UB_k$.
It is then a natural question how to prescribe the value for the offset $\tau$.
Although a general rule seems out of reach, it is always possible---if the incumbent solution is not optimal---to find a value $\tau>0$ such that the feasible set $\XX$ intersected with conditions \eqref{eq:offset_constraints} is nonempty.
Since the offset $\tau$ has to be sufficiently small for feasibility but large to promote fast decrease for the objective,
we adopt a backtracking procedure to adaptively select $\tau_k$, in analogy with the PGM outlined in \cref{alg:PGM}.
This feature is enabled in \cref{alg:BUI} with the flag \texttt{use\_offset}.

\subsection{Cuts at lower bound}\label{sec:cut_lower_bound}

When a local solver is invoked, it is possible that $x^k_{\LB} \ne x^{k+1}$.
In this case, a tighter relaxation can be obtained by adding a cut at both points $x^{k+1}$ and $x^k_{\LB}$.
However, preliminary numerical results have shown that \emph{always} adding both cuts does not improve performance;
the benefit of a possibly tighter relaxation is surpassed by the increased cost to solve subproblems with more constraints.
As an intermediate strategy, we add the cut at the lower bound point $x^k_{\LB}$ only if the condition
\begin{equation}\label{eq:LB_cuts_descent}
	\innprod{\nabla f(x^{k+1})}{x^k_{\LB}-x^{k+1}}
    \leq
    0
\end{equation}
is satisfied.
Condition \eqref{eq:LB_cuts_descent} gives an indication of whether the two points $x^k_{\LB}$ and $x^{k+1}$ are on opposite sides with respect to optimal solutions or not; see \cref{fig:cvx_illustration}.
If they are, adding both cuts can deliver a tighter relaxation;
if they are not, only the cut at $x^{k+1}$ is added
(having a lower objective value than $x^k_{\LB}$, it likely yields a stronger cut).
This feature is enabled in \cref{alg:BUI} with the flag \texttt{use\_LB\_cuts}.

\section{Numerical experiments}\label{sec:numerical_experiments}

We now present numerical results for variants of \cref{alg:BUI}.
These procedures were implemented in \texttt{python} using the open-source package \texttt{HiGHS} as MILP solver.\footnote{All tests were conducted on a laptop with a 2.80 GHz Intel Core i7 processor, 16 GB RAM, running \texttt{python} 3.11.5 and \texttt{highspy} 1.10.0.}
The baseline scheme is without local solver (CPM), otherwise we use the projected gradient scheme (PGM) as local solver.
In the latter case, it is also possible to offset the cuts and to include a cut at the lower bound, as discussed respectively in \cref{sec:cut_offset,sec:cut_lower_bound}.
Overall, these options result in five distinct solver configurations, see \cref{tab:solver_config},
which share the following setup:
tolerance $\varepsilon=10^{-9}$, initial offset $\tau_0=\infty$, gap factor $\kappa_g = 0.1$, backtracking factor $\kappa_\tau = 0.5$,
time limit at 100 seconds.
All solver parameters for the mixed-integer linear programs were set to their defaults.
For \cref{alg:PGM} we set $\alpha=10^{-3}$, $\beta=\nicefrac{1}{2}$, $\gamma_0=1$ and check only \eqref{eq:PGqualityratio} at \stepref{step:PGM:linesearch:if}.

\begin{table}
    \renewcommand{\arraystretch}{1.1}
	\centering%
	\begin{tabular}{c|ccc}%
		& ~\texttt{use\_local\_solver}~ & ~\texttt{use\_offset}~ & ~\texttt{use\_LB\_cuts}~ \\
		\hline
		CPM & $\cdot$ & $\cdot$ & $\cdot$ \\
        \hline
		PGM & \checkmark & $\cdot$ & $\cdot$ \\
		PGM-$\tau$ & \checkmark & \checkmark & $\cdot$ \\
		PGM-LB & \checkmark & $\cdot$ & \checkmark \\
		PGM-$\tau$-LB & \checkmark & \checkmark & \checkmark \\
        \hline
	\end{tabular}
	\caption{Solver configurations for \cref{alg:BUI}. The symbol \checkmark indicates that the corresponding flag is set to \texttt{true}; otherwise it is \texttt{false}. \cref{alg:PGM} acts as the local solver (PGM).}%
	\label{tab:solver_config}%
\end{table}

We will compare solver configurations on benchmark instances of binary quadratic problems with cardinality constraints, particularly maximum diversity problems (MDPs), from \cite{marti2010branch,lima2017solution}.\footnote{Available online at \url{https://grafo.etsii.urjc.es/optsicom/mdp.html} and \url{https://sites.google.com/site/cbqppaper/}.}
We refer to \cite{bruglieri2006annotated} for an overview of combinatorial problems subject to cardinality constraints and their applications.

\subsection{Metrics and profiles}

The performance of various solver configurations is assessed on some sets of test instances and based on how quickly the best objective values are reached.
Since \cref{alg:BUI} generates only feasible iterates, we can track the best objective $\UB_k$ along the iterations and compare how much and how quickly it decreases with each solver.
Specifically, let $P$ and $S$ denote a set of problems and a set of solvers, respectively.
For each problem $p\in P$ and solver $s\in S$, we record the best value $f_{s,p}$ and the wall-clock runtime $t$ in each iteration $k$.
For each problem $p$, we store the (common) initial value $f_p^0 \coloneqq f_{s,p}(0)$ and estimate the best known value $f_p^\star$---from the literature, if available, otherwise as $\min_{s,t} f_{s,p}(t)$.
Then, for each solver $s$ and problem $p$ we can define the residue function $\resid_{s,p}\colon [0,\infty)\to[0,1]$, mapping a (iteration or runtime) budget $t\geq 0$ to 
\[
\resid_{s,p}(t)
\coloneqq
\frac{f_{s,p}(t)-f_p^\star}{f_p^0-f_p^\star}
.
\]
By construction, the residue $\resid_{s,p}$ satisfies $\resid_{s,p}(0)=1$, $\resid_{s,p}(t)\in [0,1]$ for all $t\geq 0$, and it monotonically decreases with $\lim_{t\to\infty} \resid_{s,p}(t)\geq 0$.
Then, more efficient algorithms are those for which the residue decreases the most (reaching the best objective value) and faster (that is, with less budget).

Using the residue defined above, each pair $(p, s)$ corresponds to a residue profile, namely the graph of $\resid_{s,p}$.
Then, the behavior of a solver over the whole set of problems $P$ can be assessed based on the statistical properties of all its residue profiles.
The main question is:
given a computational budget, how well is any problem going to be solved?
We address this issue with the representation of \emph{median residue profiles}, which capture the distribution of the residual profiles for each solver $s\in S$.
These show, at any given budget $t\geq 0$, the median and the interquartile range of the set $\{ \resid_{s,p}(t) \,|\, p\in P \}$.

Another statistical description is given by the \emph{residue distribution} at some given budget $t\geq 0$, which depicts the cumulative distribution of residue across all problems in the test set.

\subsection{Maximum diversity problems}

The maximum diversity problem (MDP) consists of selecting a subset of $m$ elements from a set of $n$ elements in such a way that the sum of the distances between the chosen elements is maximized.
The definition of distance between elements is customized to specific applications.
Consider \eqref{eq:P} with quadratic $f$ and set $\XX$ binary and with a cardinality constraint:
\begin{equation}\label{eq:cbqp}
	\minimize_x~ \frac{1}{2} x^\top Q x
	\qquad
	\stt~x \in \{0,1\}^n ,\quad \sum_{i=1}^n x_i = m ,
\end{equation}
where the symmetric matrix $Q\in\R^{n\times n}$ and the integer $m \in (0,n)$ are from the publicly available \texttt{MDPLIB} benchmark library \cite{marti2010branch}.
Following \cite{bui2024solving}, we use the \texttt{GKD} data set, which consists of 145 matrices whose entries were calculated as Euclidean distances from randomly generated points with coordinates in the range of 0 to 10.
This includes three subsets of instances (\texttt{GKD-a}, \texttt{GKD-b}, \texttt{GKD-c}), whose size spans from $n=10$ to $n=500$, and cardinality between $m=2$ and $m=50$.
As common initial guess for all solvers we make $x^0$ feasible with the first $m$ entries equal to one.

\begin{figure}
	\centering%
    \includegraphics{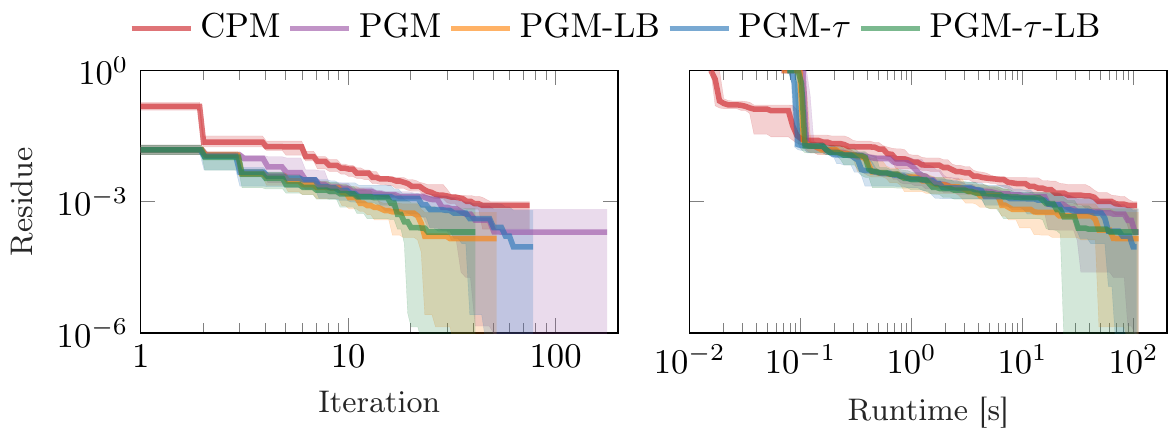}\\
	\includegraphics{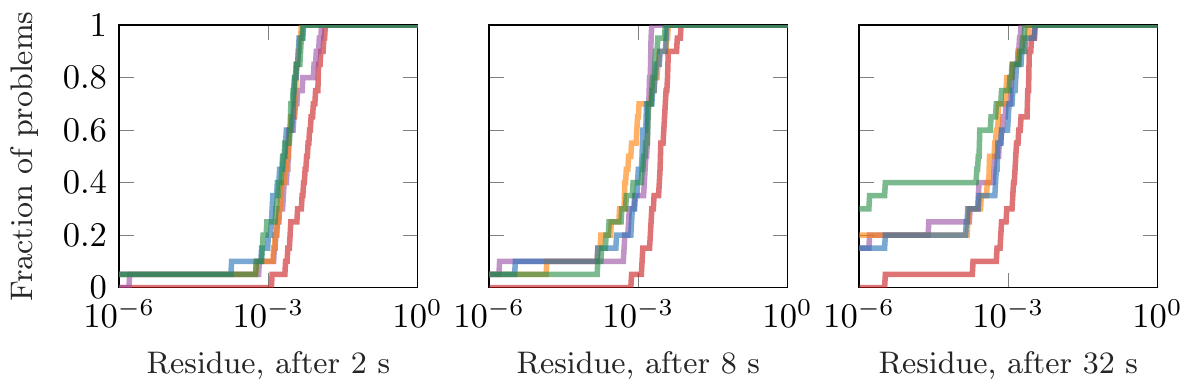}%
	\caption{Comparison of solver configurations on \texttt{GKD-c} instances: median residue profiles relative to iterations (top left panel) and runtime (top right panel), and residue distributions at fixed runtimes (bottom panels). Residue calculation based on the best known values from \texttt{MDPLIB}.}%
	\label{fig:GKD_c}%
\end{figure}

\paragraph{Results}
We begin with the largest instances, \texttt{GKD-c}, for which the best known values are available from \texttt{MDPLIB}.
Median residue profiles and residue distributions are depicted in \cref{fig:GKD_c}.
Analogous results have been obtained on test sets \texttt{GKD-a} and \texttt{GKD-b}; see \cref{sec:additional_num_results}.
These figures show that all PGM variants tend to find better values and with fewer iterations.
Indeed, neglecting an initial overhead, all PGM variants require less runtime, despite the additional effort needed to solve the subproblems.
The two variants using LB cuts seem also slightly faster than those without, but this edge might not be statistically significant.
The plain PGM configuration tends to be slower and less effective than those variants with additional tightenings.

Another observation is that the PGM variants have a wider distribution of residues toward the smaller values; see the top panels in \cref{fig:GKD_c} for instance.
This phenomenon suggests that, at least for some instances, the PGM variants are very effective heuristics, significantly accelerating the attainment of lower objective values.
The shaded bands in the median residue profiles are not standard error bars, but rather interquartile range over all instances, so the difference between the solvers' lines is indeed significant.
Analogous conclusions are valid for the \texttt{GKD-a} and \texttt{GKD-b} test sets, although with weaker statistical significance.

\subsection{Lima--Grossmann problems}
A test set for boolean quadratic programming (BQP) with cardinality constraint was generated also by Lima and Grossmann \cite{lima2017solution}.
Although the provided matrices yield possibly nonconvex functions $f$ in \eqref{eq:cbqp},
we apply \cref{alg:BUI} with the sole precaution of using cuts regularized as in \cite{spiers2022exact} (exploiting the conditional convexity of $f$ over $\XX$).
The test instances comprise five families of matrices $Q$ and five matrices are generated for each family, resulting in 80 instances.
Following \cite{lima2017solution}, we consider the cases $(n,m)=(50,10)$ and $(n,m)=(50,40)$.

\begin{figure}[tbh]
	\centering%
	\includegraphics{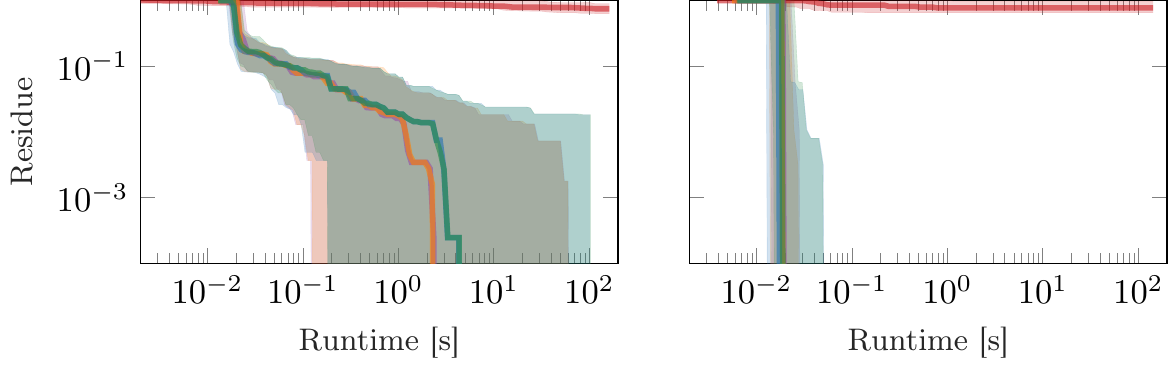}%
	\caption{Comparison of solver configurations on Lima--Grossmann problems: median residue profiles relative to runtime. Test instances with $n=50$ and $m=10$ (left) and $m=40$ (right). Legend as in \cref{fig:GKD_c}.}%
	\label{fig:LimaGrossmann}%
\end{figure}

\paragraph{Results}
The median residue profiles for both cases are depicted in \cref{fig:LimaGrossmann}: it appears that CPM makes little progress compared to the PGM variants,
which can vastly benefit from exploiting also the nonconvex function $f$ to heuristically deliver local improvements.
The residue distributions reported in \cref{fig:LG_distributions} give a complementary account of this significant advantage of PGM variants over CPM.

\section{Conclusions}
This paper presents a new approach for strengthening cuts in cutting-plane methods.
At each iteration, we use a local solver to identify local solutions to generate stronger cuts, while preserving the direction obtained from the cutting plane model.
Determining how many cuts to add and where to place them remains critical for overall performance.
Future research will focus on quantifying the strength of a cutting plane oracle rather than evaluating individual cuts in isolation, as oracle level measures can provide deeper insight into the performance and applicability of cutting plane methods across a broad class of mixed-integer nonlinear problems.

\subsubsection*{Acknowledgements}
\TheAcknowledgements

\bibliographystyle{habbrv}
\phantomsection%
\addcontentsline{toc}{section}{References}%
\bibliography{references}

\appendix

\section{Stationarity and criticality}\label{sec:local_optimality}

Denoting $\limnormalcone_A(\bar{x})$ the limiting normal cone to $A$ at $\bar{x}$, a point $\bar{x}\in\XX$ is called \emph{stationary} for \eqref{eq:localProblem} if
\begin{equation}\label{eq:stationarity}
	- \nabla f(\bar{x}) \in \limnormalcone_A(\bar{x}) .
\end{equation}
Notice that stationarity requires feasibility, namely $\bar{x} \in A$, for $\limnormalcone_A(\bar{x})$ to be nonempty.
This stationarity condition is not practical though.
It is easier to check the criticality condition given in \eqref{eq:criticality},
which is also stronger since, for all $v, z \in \R^n$,
\begin{equation}
	z\in\projection_A(v)
	\quad\Rightarrow\quad
	v-z\in\limnormalcone_A(z) .
\end{equation}
The converse relation holds true if $A$ is convex, which is not the case here.
See \cite[Prop. 3.5]{themelis2018forward} for a characterization of stationarity, criticality and optimality.

\paragraph{Proof of \cref{prop:PGM:finiteExactTermination}}
\begin{proof}
	Interpreting the Euclidean projection onto $A$ as the proximal mapping of the indicator $\indicator_A$ of $A$, \cref{alg:PGM} is an instance of a monotone, adaptive proximal gradient method,
    a special case of \cite[Alg. 3.1]{demarchi2023monotony}.
    Then, the iterations are well-defined and, up to neglecting the first iteration, it is $x^j \in A$ and $f(x^{j+1}) \leq f(x^j)$ for all $j\in\N$.
	Moreover, it follows from \cite[Lemma 4.3(iii)]{demarchi2023monotony} that $\| x^{j+1} - x^j \| \to 0$, regardless of the initial point $x^0$.
	Therefore, since $A\subset \N^n$, it must be $x^{J+1} = x^J$ for some finite $J\in\N$.
	Combining with \stepref{step:PGM:projection} of \cref{alg:PGM}, it is $x^J = x^{J+1} \in \projection_A( x^J - \gamma_{J+1} \nabla f(x^J) )$, corresponding to the criticality condition \eqref{eq:criticality} for \eqref{eq:localProblem}, with $\eta = \gamma_{J+1} > 0$.
    \qed
\end{proof}

\section{Trust-region methods}\label{sec:trust_region}

The projection-based scheme presented in \cref{sec:projected_gradient} is practical as long as the projection operator can be easily evaluated.
When the set $A$ is binary, projecting onto $A$ corresponds to solving a MILP,
but more general sets require solving a mixed-integer \emph{quadratic} program (MIQP), which can become prohibitive in practice.
Then, for \eqref{eq:localProblem} one can look at a linear model of the objective subject to $x\in A$ as well as a trust-region (TR) constraint $\|x-x^j\|_p \leq \Delta_j$.
The latter forces updates to remain in a region around the current estimate $x^j$,
controlled by the trust-region radius $\Delta_j$.
The PGM iteration \eqref{eq:PGiteration} can then be replaced with the TR update
\begin{equation}\label{eq:TRiteration}
	x^{j+1} \in \arg\min_{x\in A} \left\{
    \innprod{\nabla f(x^j)}{x}
    \,\middle\vert\,
    \|x-x^j\|_p \leq\Delta_j
    \right\},
    \quad j=0,1,\ldots
\end{equation}
with some suitable radius $\Delta_j > 0$.
When set $A$ is described by linear inequalities,
selecting a polyhedral norm (i.e. $p\in\{1,\infty\}$) results in \eqref{eq:TRiteration} being a MILP.
As in \cref{alg:PGM}, a backtracking procedure can be adopted to find a radius $\Delta_j>0$ that delivers sufficient decrease, in the sense of \eqref{eq:PGsuffdecrease} or \eqref{eq:PGqualityratio}.
See \cite{demarchi2025mixed} for a TR approach to deal with mixed-integer variables as well as for a detailed algorithm and convergence analysis.

\section{Additional numerical results}\label{sec:additional_num_results}

\begin{figure}[tbh]
	\centering%
	\includegraphics{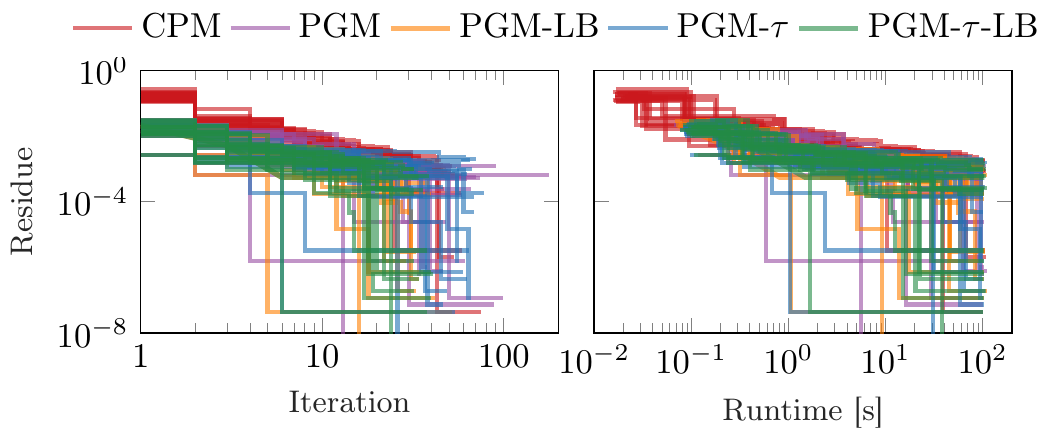}%
	\caption{Comparison of solver configurations on \texttt{GKD-c}: residue along iterations (left) and runtime (right) for each problem instance.}%
	\label{fig:GKD_c_logger}%
\end{figure}

\begin{figure}[tbh]
	\centering%
    \includegraphics{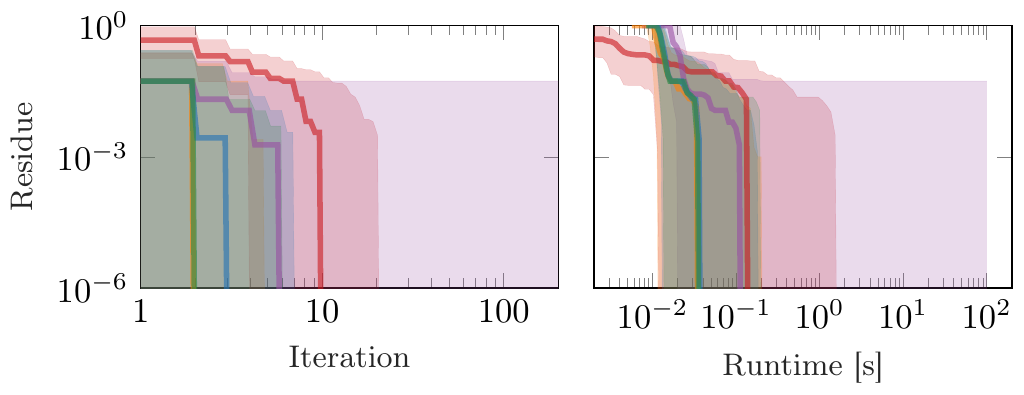}\\
	\includegraphics{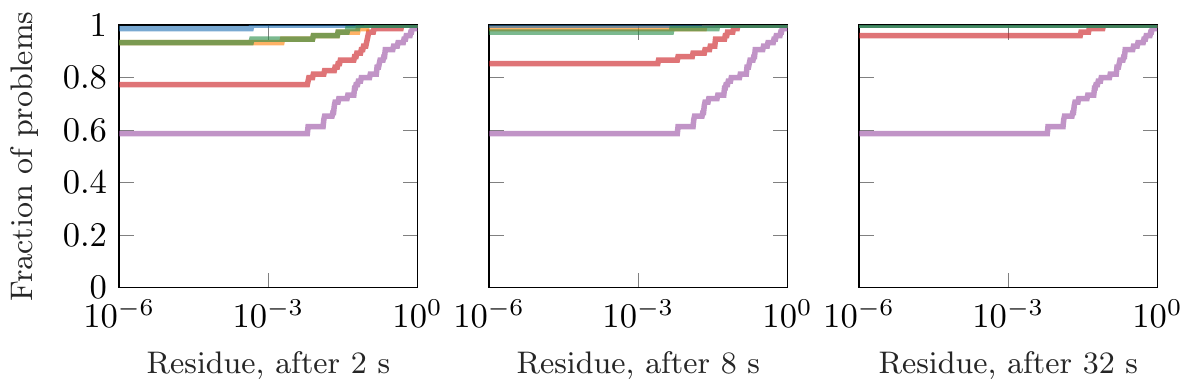}%
	\caption{Comparison of solver configurations on \texttt{GKD-a} instances: median residue profiles relative to iterations (top left panel) and runtime (top right panel), and residue distributions at fixed runtimes (bottom panels). Legend as in \cref{fig:GKD_c_logger}.}%
	\label{fig:GKD_a}%
\end{figure}

\begin{figure}[tbh]
	\centering%
    \includegraphics{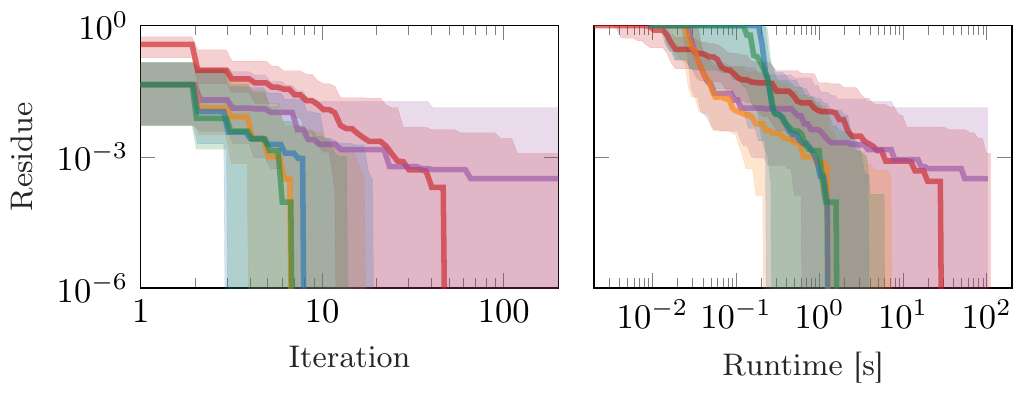}\\
	\includegraphics{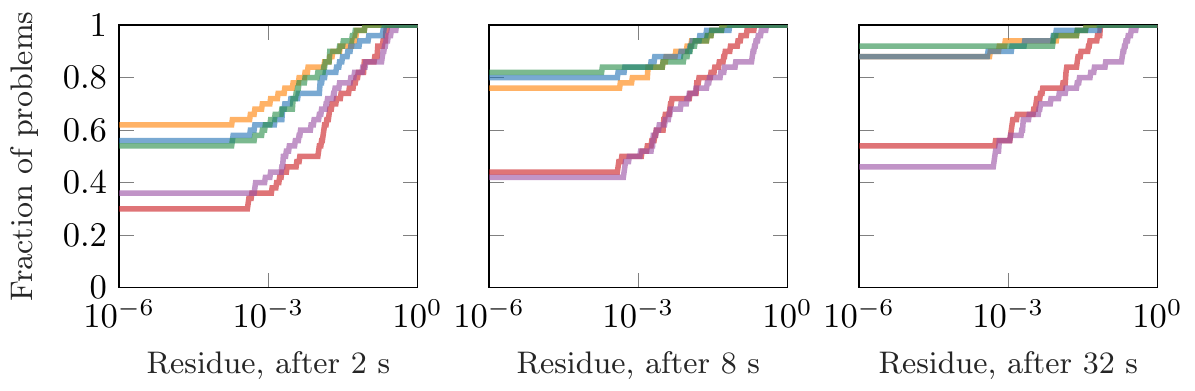}%
	\caption{Comparison of solver configurations on \texttt{GKD-b} instances: median residue profiles relative to iterations (top left panel) and runtime (top right panel), and residue distributions at fixed runtimes (bottom panels). Legend as in \cref{fig:GKD_c_logger}.}%
	\label{fig:GKD_b}%
\end{figure}

An exemplary illustration of residue profiles $\resid_{s,p}$ as a function of computational budget is depicted in \cref{fig:GKD_c_logger};
these are graphically summarized by the median residue profiles in \cref{fig:GKD_c}.
An analogous statistical representation for the test sets \texttt{GKD-a} and \texttt{GKD-b} is given in \cref{fig:GKD_a,fig:GKD_b},
from which similar conclusions can be drawn, despite the larger variability.

\begin{figure}[tbh]
	\centering%
	\includegraphics{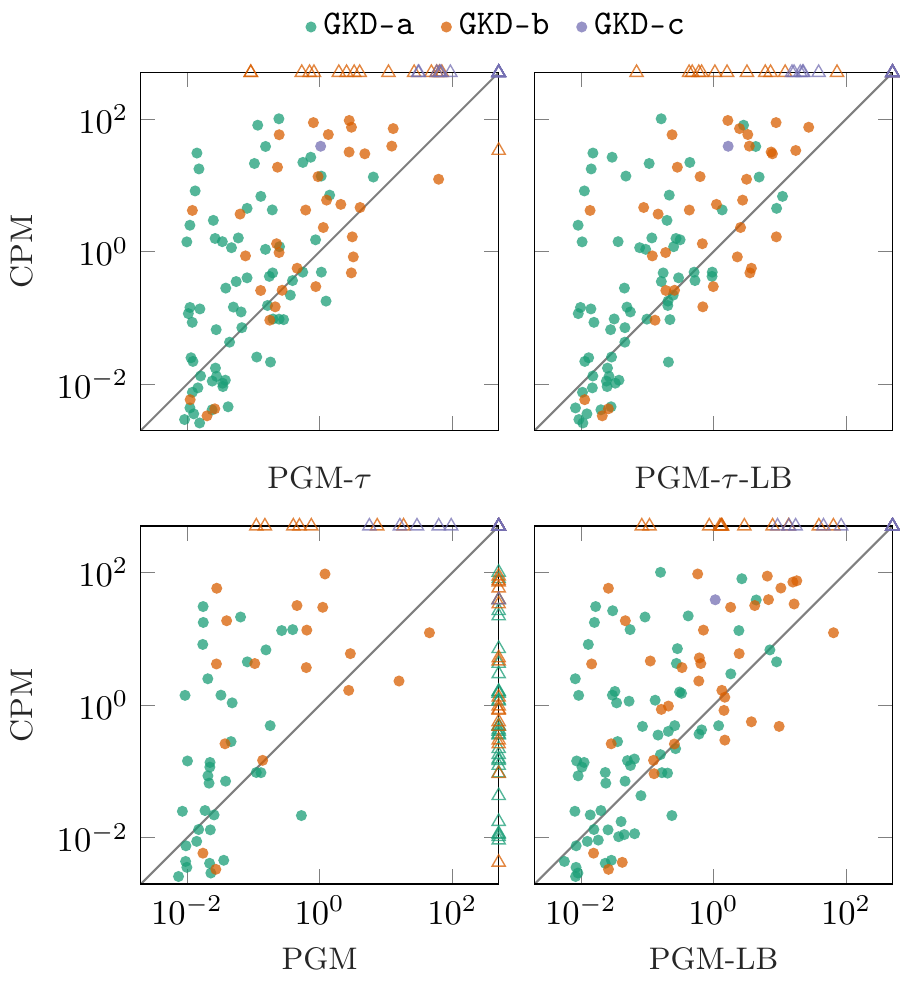}%
	\caption{Comparison of CPM against PGM variants on \texttt{MDPLIB} instances: runtime to reach the target residue $\resid^* = 10^{-6}$. If either solver fails, the pair of runtimes is marked with a triangle (with a large runtime for the failing solver).}%
	\label{fig:tbest_cpm}%
\end{figure}

\begin{figure}[tbh]
	\centering%
	\includegraphics{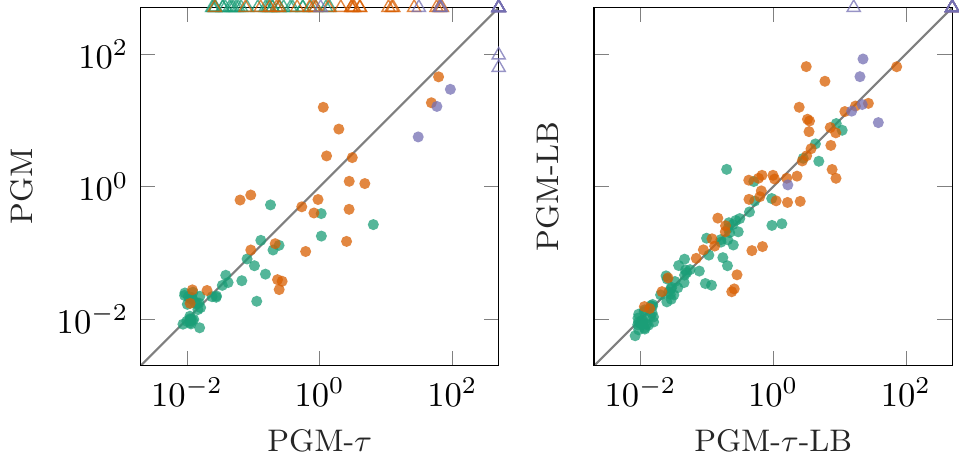}%
	\caption{Comparison of PGM variants with and without tightening offset $\tau$ on \texttt{MDPLIB} instances: runtime to reach the target residue $\resid^* = 10^{-6}$. If either solver fails, the pair of runtimes is marked with a triangle (with a large runtime for the failing solver). Legend as in \cref{fig:tbest_cpm}.}%
	\label{fig:tbest_effect_offset}%
\end{figure}

\begin{figure}[tbh]
	\centering%
	\includegraphics{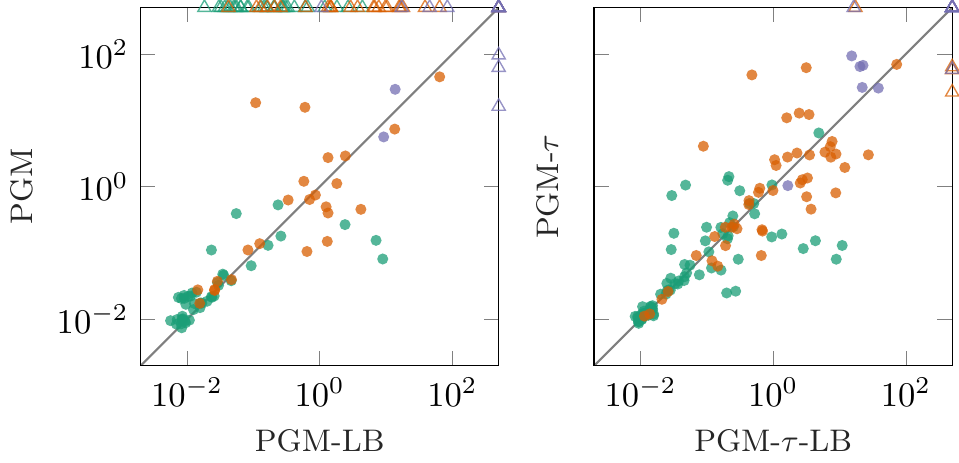}%
	\caption{Comparison of PGM variants with and without tightening LB cuts on \texttt{MDPLIB} instances: runtime to reach the target residue $\resid^* = 10^{-6}$. If either solver fails, the pair of runtimes is marked with a triangle (with a large runtime for the failing solver). Legend as in \cref{fig:tbest_cpm}.}%
	\label{fig:tbest_effect_LB}%
\end{figure}

A close comparison of the original CPM against the PGM variants is reported in \cref{fig:tbest_cpm} in terms of runtime to reach a target residue $\resid^* = 10^{-6}$.
The scattered points are, for the most part, above the diagonal, which means that CPM often takes longer than all PGM variants to decrease $\UB_k$ and reach the target residue.
This observation is valid across budget metrics (iterations and runtime) and test sets (\texttt{GKD-a}, \texttt{GKD-b} and \texttt{GKD-c}).
However, from the comparison of CPM against PGM, it appears that the latter does not perform particularly well either, as it often fails to reach the desired residue.
In contrast, PGM variants with tightenings are much more effective: almost all test instances are solved up to the target residue within the time limit.
This assessment also reflects the residue distributions in \cref{fig:GKD_a,fig:GKD_b,fig:GKD_c}.
However, it seems that adopting both tightening techniques does not yield further improvements.
Pairwise comparisons of PGM variants to highlight the marginal effects of offset and LB cuts are depicted in \cref{fig:tbest_effect_offset,fig:tbest_effect_LB} respectively.
In both cases, the beneficial effect is clear and significant over the plain PGM (left panels), but not when the other tightening procedure is already enabled (right panels).

\begin{figure}[tbh]
	\centering%
    \includegraphics{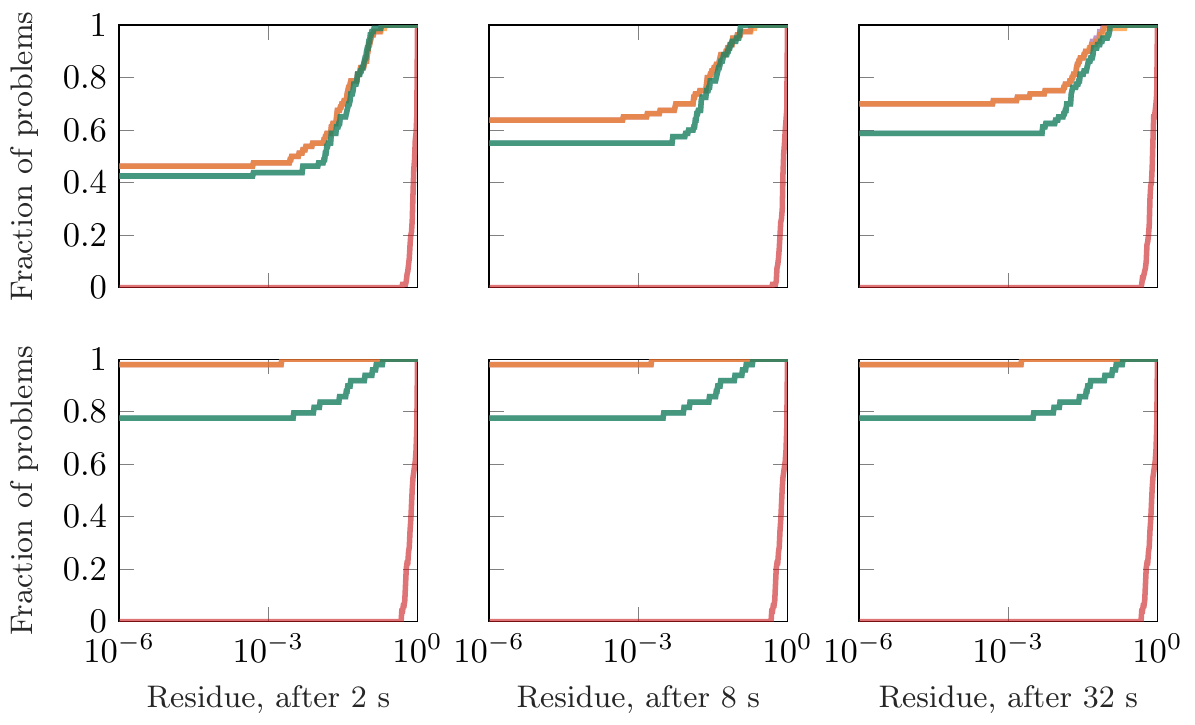}%
	\caption{Comparison of solver configurations on Lima--Grossmann problems: residue distributions relative to runtime. Test instances with $(n,m)=(50,10)$ (top panels) and $(n,m)=(50,40)$ (bottom panels). Note that curves for solver configurations with and without LB cuts overlap. Legend as in \cref{fig:GKD_c_logger}.}%
	\label{fig:LG_distributions}%
\end{figure}

Finally, we report the residue distributions obtained for the Lima-Grossmann instances in \cref{fig:LG_distributions}, corresponding to the median residue profiles of \cref{fig:LimaGrossmann}.
The main observation is the vast improvement over the CPM solver, but it is interesting to notice also the following.
\begin{itemize}
    \item Distributions of solver configurations with and without LB cuts overlap, meaning that \texttt{use\_LB\_cuts} does not affect the performance in this case.
    This happens because condition \eqref{eq:LB_cuts_descent} is never satisfied, hence no LB cuts are added.
    \item PGM variants with enabled \texttt{use\_offset} are slower than those without, indicating that the offset $\tau$ does not contribute significantly to tighten the local subproblem in this case.
\end{itemize}

\end{document}